\mag=\magstep1
\documentstyle{amsppt}\input amsppt1
\pageheight{24true cm}
\pagewidth{16true cm}

\parindent=4mm\parskip=3pt plus1pt minus.5pt 
\nologo\NoRunningHeads\NoBlackBoxes

\def\e{\hookrightarrow}

\def\f{\flushpar }
\def\nl{\newline }
\def\np{\newpage }
\def\x{\times }

\def\o{ordinary sense slice 1-link}
\def\oo{ordinary sense slice 1-link }

\topmatter
\title   
Some properties of ordinary sense slice 1-links:
Some answers to the problem (26)  of Fox  
\endtitle
\author
Eiji Ogasa$^*$
\endauthor
\thanks{
{\it 1991 Mathematics Subject Classification.} Primary 57M25, 57Q45 \nl
 Keywords. ordinary sense slice 1-links, the Arf invariants, 
$n$-dimensional knots and links, Suzuki-Terasaka  diagrams, 
realizable 4-tuple of  links \nl
$^*$This research was partially supported by Research Fellowships 
of the Promotion of Science for Young Scientists.  \newline
This paper is published in \newline
Proceedings of the American Mathematical society  126, 1998, P.2175-2182.\newline
This manuscript is not the published version.
}
\affil
Department of Mathematical Sciences,  University of Tokyo\\
  Komaba, Tokyo 153,    Japan\\
ogasa\@ms.u-tokyo.ac.jp, ogasa\@ms513red.ms.u-tokyo.ac.jp
\endaffil
\endtopmatter

\document

\baselineskip11pt
{\bf Abstract.}
We prove that, 
for any ordinary sense slice 1-link $L$, we can define the Arf invariant 
and Arf($L$)=0. 
We prove that,  
for any $m$-component 1-link $L_1$, there exists a $3m$-component \oo $L_2$ 
of which $L_1$ is a sublink.

\head 1. Introduction and Main results \endhead

In  
[3] %[Fox]
Fox submitted  the following problem about 1-links.  
Here, note that ``slice link''  in the following problem is now  called 
``ordinary sense slice link,''   
and  
``slice link in the strong sense'' in the following problem is now called 
``slice link'' by knot theorists.

\definition{Problem 26 of 
[3] %[Fox]
} Find a necessary condition for $L$ 
to be a slice link; a slice link in the strong sense.  
\enddefinition

Our purpose is to give some answers to the former part of this problem.
The latter half is not discussed here. 
The latter half seems discussed much more often than the former half. 
See e.g. 
[2], %[Cochran and Orr]
[5] %[Gilmer and Livingston]
 and 
[14]. %[Levine 1994]

We review the definition of ordinary sense slice links   
 and that of slice links, which we now use.

 We suppose $m$-component 1-links are oriented and ordered. 

Let $L=(K_1,...,K_m)$ be a $m$-component 1-link in $S^3=\partial B^4$. 
$L$ is called a {\it slice 1-link}, 
which is ``a slice link in the strong sense'' in the sense of Fox, 
if there exist 2-discs $D^2_i (i=1,...,m)$ in $B^4$ such that 
$D^2_i\cap\partial B^4$ =$\partial D^2_i$, $D^2_i\cap D^2_j$=$\phi(i\neq j)$, 
and $(\partial D^2_1,...,\partial D^2_m)$ in $\partial B^4$ defines $L$.

Take a 1-link $L$ in $S^3$. 
Take $S^4$ and regard $S^4$ as $(\Bbb R^3\x \Bbb R)\cup \{\infty\}$. 
Regard the 3-sphere $S^3$ as $\Bbb R^3\cup \{\infty\}$ in $S^4$. 
$L$ is called an {\it \o}, 
which is ``a slice link'' in the sense of Fox, 
 if there exists an embedding 
$f:S^2\e \Bbb R^3\x \Bbb R$ such that $f$ is transverse to 
$\Bbb R^3\x \{0\}$ and 
$f(S^2)\cap(\Bbb R^3\x \{0\})$ in $\Bbb R^3\x \{0\}$ defines $L$.  
Suppose $f$ defines a 2-knot $X$. 
Then $L$ is called a {\it cross-section} of the 2-knot $X$.

From now on we use the terms in the sense of the present.

Ordinary sense slice 1-links have the following properties. 
\proclaim{Theorem 1.1}
Let $L$ be an \o.
Then the followings hold. 
\roster
\item 
$L$ is a proper link.
\item
Arf($L$)=0. 
\endroster
\endproclaim

\definition {Note}
Although our proof proves (1) and (2) simultaneously, 
once (1) is known (2) follows easily from the known result 
that a proper link which is an ordinary sense slice link 
has trivial Arf invariant.  
See  e.g.  
[4], [9], [16] and [19]. %[Kawauchi][Murasugi][Robertello][Gilmer].  
\enddefinition

\proclaim{Theorem 1.2}
For any $m$-component 1-link $L$, there exists a $3m$-component \oo $L'$ 
of which $L$ is a sublink. 
\endproclaim

\definition {Note}
When $L$ is not slice and $m$=1, it is obvious that `$3m$' is best possible. 
When $m\geqq2$,  `$3m$' is not best possible 
even if no components of $L$ are slice knots. 
See the example for Note 3.3 and Figure IV  in \S 3. 
\enddefinition 

 Theorem 1.2 follows from Theorem 1.3 obviously.
\proclaim{Theorem 1.3}
For any $m$-component 1-link $L$=$(K_1,...,K_m)$, 
there exists an embedding 
$f:S^2_1\amalg...\amalg S^2_m\e \Bbb R^3\x \Bbb R$ and 
a $3m$-component 1-link $L'$ with the following properties. 
 \roster
 \item
  $f$ is transverse to $\Bbb R^3\x \{0\}$ and 
$f(S^2_1\amalg...\amalg S^2_m)\cap (\Bbb R^3\x \{0\})$ in $\Bbb R^3\x \{0\}$ 
defines $L'$. 
  \item
  $L$ is a sublink of $L'$. 
  \item 
  $K_i\subset f(S^2_i) (i=1,...,m)$. 
 \endroster

\endproclaim

This paper is organized as follows. 
In \S 2 we prove Theorem 1.1 by using a result of the author's 
[17]%[OgasaI]. 
In \S 3 we review Suzuki-Terasaka diagrams of 1-links 
and the fact that any 1-link has a  Suzuki-Terasaka diagram (Theorem 3.1).
We use this diagram to show Theorem 1.3 and Theorem 3.2.

\head  
2. The proof of Theorem 1.1
\endhead
In order to  prove Theorem 1.1, we review a result of the author in 
[17]. %[OgasaI]. 

\f {\bf Definition }
$T=(L_1, L_2, X_1, X_2)$ is called a 
{\it 4-tuple of links } if the following conditions (1), (2) and (3) hold.

\f (1) $L_i=(K_{i1},...,K_{im_i})$  is an oriented ordered $m_i$-component 
1-dimensional link  $(i=1,2).$  \newline

\f (2)    $m_1=m_2.$

\f (3)     $X_i$ is an oriented 3-knot. 

\f {\bf Definition }
A 4-tuple of links $(L_1, L_2, X_1, X_2)$ 
is said to be  {\it realizable } if there exists 
a smooth transverse immersion  
$f:S^3_1\coprod S^3_2 \looparrowright S^5$  
satisfying the following conditions (1) and (2). 

\f(1)   $f\vert S^3_i$ is a smooth embedding 
and defines the  3-knot $X_i (i=1,2)$ in $S^5$.  

\f(2)   For  $C=f(S^3_1)\cap f(S^3_2)$, 
the inverse image  $f^{-1}(C)$ in $S^3_i$ defines the  1-link 
$L_i (i=1,2).$ Here, the orientation of $C$ is induced 
naturally from the preferred orientations of $S^3_1, S^3_2,$ and $S^5$,   
 and an arbitrary order is given to the components of $C$.

The following theorem characterizes the realizable 4-tuples of links.

\proclaim{Theorem 2.1}
A 4-tuple of links $T=(L_1, L_2, X_1, X_2)$ is realizable  if and only if 
 $T$  satisfies one of the following conditions i) and  ii). 

i) Both $L_1$ and $L_2$ are proper links, and 

$$\roman{Arf}(L_1) = \roman{Arf}(L_2).$$ 

ii) Neither $L_1$ nor $L_2$ 
is proper, and 
$$ lk(K_{1j}, L_1-K_{1j})
\equiv
  lk(K_{2j}, L_2-K_{2j}) \quad \roman{mod\hskip2mm  2}
 \quad \roman{for \hskip2mm all}\hskip2mm j. $$ 

\endproclaim

Note. In 
[18], %[Ogasa II], 
the author discussed high dimensional version of Theorem 2.1.

We begin the proof of Theorem 1.1.

Take 
 $f:S^2\e \Bbb R^3\x \Bbb R$  
 such that  $L$ is  a cross section of the 2-knot defined by $f$.

Regard $\Bbb R^3\x \Bbb R$ as  
$\Bbb R^3\x \Bbb R\x \{0\}$  $\subset$ $\Bbb R^3\x \Bbb R\x \Bbb R$. 
Make $S^5$ 
from  $\Bbb R^3\x \Bbb R\x \Bbb R$  
by the one point compactification.
Here, the 4-sphere $S^4$ is 
$(\Bbb R^3\x \Bbb R\x \{0\}\cup \{\infty\})$  $\subset$ 
$S^5$=$(\Bbb R^3\x \Bbb R\x \Bbb R\cup \{\infty\})$. 
There exists a 3-knot $X_1$ in $\Bbb R^3\x \Bbb R\x \Bbb R$
such that 
$X_1$ $\cap$ $\Bbb R^3\x \Bbb R\x \{0\}$ coincides with $f(S^2)$ 
because all 2-knots are slice by a theorem of Kervaire in 
[10]. %[Kervaire]. 
$(\Bbb R^3\x \{0\}\x \{0\})\cup \{\infty\}$ in $S^5$ is called a 3-knot $X_2$.
An immersion 
$g:S^3_1\coprod S^3_2 \looparrowright S^5$  such that 
$g(S^3_i)$ coincides with the above $X_i$ ($i=1,2$) 
realizes a pair of 1-links 
($X_1\cap X_2$ in $X_1$, $X_1\cap X_2$ in $X_2$). 
Here, it is obvious that 
$X_1\cap X_2$ in $X_1$is the trivial 1-link and 
 $X_1\cap X_2$ in $X_2$ is $L$. 
Therefore, by Theorem 2.1, 
$L$ is a proper link and Arf($L$)=0.

\head 3.  The proof of Theorem 1.3 \endhead

We first review the following canonical diagrams of 1-links.   

Take a $m$-component 1-link $L$. 
The set $X$ in 
$\Bbb R^3$=$\{(x, y, z)\}$ ($\subset S^3$) 
defining $L$ is called a 
{\it Suzuki-Terasaka (canonical) diagram of $L$} 
if $X$ is made as follows. 
Let $Y_i$ ($i=1,...,m$)  be the boundary of 
$\{(x, y, z)\vert$ $ i\leqq x\leqq (i+0.9), 0\leqq y\leqq 1,  z=0\}$.  
Let 
$P_{ij}$ ($j=1,...,\mu_i$) be the boundary of 
$\{(x, y, z)\vert$ $x=\frac{j}{\mu_i+1}+i, 
                      0.9\leqq y\leqq 1.1, 
                      -0.1\leqq z\leqq 0.1 \}$.  
(The orientation of $P_{ij}$ is given appropriately.  
All are not same in general.)
Let  $A_{ij}$ =
$\{(x, y, z)\vert$ $\frac{2j-1}{2\nu_i+1}+i\leqq x\leqq\frac{2j}{2\nu_i+1}+i, 
              y=0, z=0    \}$ ($j=1,...,\nu_i$).  
Let $\Sigma^m_i\mu_i$=$\Sigma^m_i\nu_i$, 
put $\lambda$ equal to this number. 
Take bands $B_l$ ($l=1,...,\lambda$) and 
make a band-sum of $Y_i$ and $P_{ij}$  by 
connecting $P_{ij}$ and $A_{i'j'}$ by $B_l$. 
(Of course $(i,j)$ does not necessarily coincides with $(i',j')$ and 
the set of $(i,j)$ coincides with the set of $(i',j')$. )
Then the band-sum is $X$.

\f{\bf Theorem 3.1.}
Any 1-link  has a Suzuki-Terasaka canonical diagram.

The proof of Theorem 3.1 is elementary.

The usefulness of the above canonical diagram 
is firstly pointed out by 
Prof. Shin'ichi Suzuki 
and Prof. Hidetaka Terasaka. 
Suzuki-Terasaka canonical diagrams are used, 
for example,  in 
[15], [23]%[Miyazaki and Yasuhara][Suzuki1969]
and 
[25]. %[Yamamoto]

In Figure I 
there is written an example of the canonical diagram.

\vskip3mm
%\hskip30mm  
Figure I. 
%See the last pages.

You can obtain this figure 
by clicking `PostScript' in the right side of 
the cite of the abstract of this paper in arXiv 
(https://arxiv.org/abs/the number of this paper). 

You can also obtain it from the author's website,  
which can be found by typing his name in search engine.

\vskip3mm

We begin the proof of Theorem 1.3.

Take a Suzuki-Terasaka canonical diagram of $L$=$(K_1,...,K_m)$  in 
$\Bbb R^3$=$\{(x, y, z)\}$ ($\subset S^3$). 
(See Figure II (1). )
Take sets $P_i$ to be 
$\{(x, y, z)\vert$ $i\leqq x\leqq(i+0.9),$ $1.05\leqq y\leqq2,$ $z=0 \}$.  
Here, we can take  $P_i$ not to  intersect with all the bands in 
the Suzuki-Terasaka canonical diagram. 
Take sets $S_i$ to be 
$\{(x, y, z)\vert$ $i\leqq x\leqq(i+0.9),$ $-1\leqq y\leqq1,$ $z=0 \}$.  
Note that $K_i$ $\cap$ $S_i$ $\neq$ $\phi$. 
(See Figure II (2). )
The arc in $K_i$ whose boundary is the points
$\{(i,0,0)\}$ and $\{(i+0.9,0,0)\}$  and which does not include 
the point $\{(i,1,0)\}$   is called $l_i$. 
Carry out the band-fusion on $S_i$ by using the band $\widetilde{B_i}$ 
whose core is $l_i$. 
Then $S_i$ splits to two pieces. 
The one including the point  $\{(i,1,0)\}$ is called $Q_i$. 
The other is called $R_i$. 
Choose the band along $l_i$ so that 
$lk(P_i, R_i)$+$lk(Q_i, R_i)$=0.  
Here $P_i$ and $S_i$ are oriented counterclockwise, 
and $Q_i$ and $R_i$ are given orientations induced from $S_i$.   
( See Figure II (3). )
Here, note that  ($Q_1,...,Q_m$) defines the 1-link $L$.   
Thus we obtain a  $3m$-component 1-link $L'$  defined by 
($P_1,...,P_m,$ $Q_1,...,Q_m$,$R_1,...,R_m$)  
such that  $L$ is a sublink of $L'$.

\vskip3mm
\hskip30mm  Figure II. %See the last pages.

You can obtain this figure 
by clicking `PostScript' in the right side of 
the cite of the abstract of this paper in arXiv 
(https://arxiv.org/abs/the number of this paper). 

You can also obtain it from the author's website,  
which can be found by typing his name in search engine.

\vskip3mm

\f{\bf Claim.}
There exists 
$f:S^2_1\amalg...\amalg S^2_m\e \Bbb R^3\x \Bbb R$ 
such that 
  $f$ is transverse to $\Bbb R^3\x \{0\}$ and 
$f(S^2_1\amalg...\amalg S^2_m)\cap (\Bbb R^3\x \{0\})$ in $\Bbb R^3\x \{0\}$ 
defines the 1-link $L'$.

Proof. 
Put an embedding 
$f:S^2_1\amalg...\amalg S^2_m\e \Bbb R^3\x \Bbb R$ 
as follows. 
$f(S^2_i)$ in $\Bbb R^3\x \Bbb R$ has two minimum-discs,  
two saddle-bands and two maximum-discs.  
The minimum-discs are 
$h^0_{i1}$=
$\{(x, y, z)\vert$ 
$i\leqq x\leqq(i+0.9),$ $-1\leqq y\leqq1, z=0   \}$ $\x\{-2\}$ 
in  $\Bbb R^3 \x \{-2\}$
and 
$h^0_{i2}$=$\{(x, y, z)\vert$ 
$i\leqq x\leqq(i+0.9),$ $1.05\leqq y\leqq2, z=0 \}$ $\x\{-2\}$  
in  $\Bbb R^3 \x \{-2\}$.  
The saddle-bands are 
 $\widetilde{B_i}$ $\x \{-1\}$ in $\Bbb R^3 \x \{-1\}$ 
and 
$\{(x, y, z)\vert$ $i\leqq x\leqq(i+0.9),$ $1\leqq y\leqq1.05, z=0 \}$ 
 $\x\{1\}$ 
 in  $\Bbb R^3\x \{1\}$.  
The maximum-discs are 
$h^2_{i1}$=$\{(x, y, z)\vert$ 
$i\leqq x\leqq(i+0.9),$ $-1\leqq y\leqq-0.1, z=0 \}$   $\x\{2\}$ 
in  $\Bbb R^3 \x \{2\}$  
and 
$h^2_{i2}$=$\{(x, y, z)\vert$ 
$i\leqq x\leqq(i+0.9),$ $0.1\leqq y\leqq2, z=0  \}$  $\x\{2\}$, 
in  $\Bbb R^3 \x \{2\}$,   
where we suppose 
$\partial \widetilde{B_i}$ $\cap$ $S_i$ =
$\{(x, y, z)\vert$ $x=i,(i+0.9),$ $-0.1\leqq y\leqq0.1, z=0 \}$. 
$f(S^2_i)$ $\cap$ $(\Bbb R^3\x \{t\})$ 
($-2<t<-1$, $-1<t<1$, $1<t<2$)  
is an ordinary cross-section.  
Then 
$f(S^2_1\amalg...\amalg S^2_m)\cap (\Bbb R^3\x \{0\})$ in $\Bbb R^3\x \{0\}$ 
defines $L'$.

Therefore the proof of Theorem 1.3 complete.

\f{\bf Note.}      
See \S 2 of 
[24] %[Suzuki1976]
for the definitions of   
`minimum-disc,' `maximum-disc,'  `saddle-band' `ordinary cross-section,' etc.

\f{\bf Note.}      
The 1-link $(Q_1, R_1)$ (or  $(P_1, Q_1)$ )
is associated with a $\theta$-graph. 
The diagram of $(Q_1, R_1)$(or  $(P_1, Q_1)$ )
is what is used in Appendix of 
[15]. %[Miyazaki and Yasuhara]
Dr. Akira Yasuhara made an alternative proof of 
[12]. %[Kinoshita and Terasaka]
and written in Appendix of 
[15]. %[Miyazaki and Yasuhara]
[11] %[Kinoshita]
is a generalization of 
[12] %[Kinoshita and Terasaka]
and 
[6].  %[Harrold and Kinoshita]

Figure III illustrates $f(S^2)$ and $L'$ 
in $\Bbb R^3\x \Bbb R$ 
in the case where $L$ is the trefoil knot. 
This method of drawing subset of $\Bbb R^3\x \Bbb R$ is often used.
See e.g.
[1], [7], [9] and [24]. 
%[Suzuki1976][Carter and Saito][Kamada][KawauchiShibuya Suzuki]

\vskip3mm
\hskip30mm  Figure III.  %See the last pages.

You can obtain this figure 
by clicking `PostScript' in the right side of 
the cite of the abstract of this paper in arXiv 
(https://arxiv.org/abs/the number of this paper). 

You can also obtain it from the author's website,  
which can be found by typing his name in search engine.

\vskip3mm

We next discuss \o s in the case when 
we restrict the 2-knots of which the \o s are cross-sections.

\proclaim{Theorem 3.2}
For any $m$-component 1-link $L=(K_1,...,K_m)$, 
there exists an embedding
$g:S^2_1\amalg...\amalg S^2_m\e \Bbb R^3\x \Bbb R$ and 
a $4m$-component 1-link $L''$ with the following properties. 
 \roster
 \item
 $g$ defines the trivial 2-link. 
 \item
  $g$ is transverse to $\Bbb R^3\x \{0\}$ and 
$g(S^2_1\amalg...\amalg S^2_m)\cap (\Bbb R^3\x \{0\})$ 
in $\Bbb R^3\x \{0\}$ defines $L''$. 
  \item
  $L$ is a sublink of $L''$. 
  \item
  $K_i\subset f(S^2_i) (i=1,...,m)$. 
 \endroster
\endproclaim

Ordinary sense slice $n$-knots ($n\geqq1$) 
 which are cross-sections  of  the trivial $(n+1)$- knots 
 are discussed in 
[13], [20], [21], and [22].
%[Levine1983][Ruberman1988][Ruberman1983][Sumners]

We begin the proof of Theorem 3.2.

In order to define $g$, take the following 2-discs $D^2_i$ 
in $\Bbb R^3\x \Bbb R$. Take $f$ in the proof of the above claim. 
$f(S^2_i)$ $\cap$ $\Bbb R^3\x (-1,2]$ are two components. 
Take the  component  of the two which includes $h^2_{i1}$, say $E_i$.
$f(S^2_i)$ $\cap$ $\Bbb R^3\x [-2,1)$ are two components. 
Take the  component  of the two which includes $h^0_{i1}$, say $E'_i$.
 Then 
 $ \overline{E_i \cup E'_i} $
  is a submanifold in $\Bbb R^3\x \Bbb R$ diffeomorphic to the 2-disc .
 We call it $D^2_i$.

 Take  $D^2_i$ $\x I$ in the tubular neighborhood of  $D^2_i$ 
 in $\Bbb R^3\x \Bbb R$. 
 Take $g$ so that 
 $g(S^2_i)$ coincides with $\partial $ ( $D^2_i$ $\x I$)  
and  
 $g$ is transverse to $\Bbb R^3\x \{0\}$. Then 
$g(S^2_1\amalg...\amalg S^2_m)\cap (\Bbb R^3\x \{0\})$ 
in $\Bbb R^3\x \{0\}$ is a $4m$-component 1-link $L''$ and 
  $L$ is a sublink of $L''$.

This completes the proof of Theorem 3.2.

Figure IV illustrates  
 $g(S^2)$ and $L''$ in $\Bbb R^3\x \Bbb R$ 
 in the case where $L$ is the trefoil knot.

\definition{Note 3.3}
We can regard Figure IV as the example we mentioned in Note under Theorem 1.2 
if we think $L$ being a union of $K_1$ and $K_2$ 
and $L''$ being a union of  $K_1$, $K_2$, $K_3$ and $K_4$, 
where $K_i$ are as in Figure IV. 
\enddefinition

\vskip3mm
\hskip30mm  Figure IV. %See the last pages.

You can obtain this figure 
by clicking `PostScript' in the right side of 
the cite of the abstract of this paper in arXiv 
(https://arxiv.org/abs/the number of this paper). 

You can also obtain it from the author's website,  
which can be found by typing his name in search engine.

\vskip3mm

Comparing Theorem 1.3 with Theorem 3.2, 
it is natural to raise the following problem. 

\definition{Problem }
Let $K$ be a non-slice knot. Does there exist a $3$-component \oo $L$ 
such that $K$ is a component of $L$ and 
$L$ is a cross-section of  the trivial 2-knot? 
\enddefinition

\f{\bf Note 3.4.}
After Dr. S. Kamada received a manuscript of this paper, 
he solved this Problem  
and 
he obtained a refined version of Theorem 1.3 and 3.2.

\np

\enddocument